\documentclass[12pt]{article}
\makeatletter \oddsidemargin  -.1in \evensidemargin -.1in
\newcommand{\singlespacing}{\let\CS=\@currsize\renewcommand{\baselinestretch}{1}\tiny\CS}
\newcommand{\oneandahalfspacing}{\let\CS=\@currsize\renewcommand{\baselinestretch}{1.25}\tiny\CS}
\newcommand{\doublespacing}{\let\CS=\@currsize\renewcommand{\baselinestretch}{1.35}\tiny\CS}

\newtheorem{rule-def}[theorem]{Rule}

\RequirePackage[dvips]{graphicx} \textheight 22.5cm
\usepackage{latexsym,epsfig,enumerate,amsmath,amsfonts,amssymb,amsbsy,amsopn,mathrsfs}
\usepackage{subfigure}
\usepackage{graphicx}

\begin{document}

\title{\bf The Non-existence of Perfect Cuboid} \author{\small
  S. Maiti$^{1,2}$ \thanks{Corresponding author, Email address: {\it
      somnath.maiti@lnmiit.ac.in/maiti0000000somnath@gmail.com
      (S. Maiti)}}\\\it $^{1}$ Department of Mathematics, The LNM
  Institute of Information Technology\\\it Jaipur 302031, India\\ \it
  $^{2}$Department of Mathematical Sciences, Indian Institute of
  Technology (BHU),\\ Varanasi-221005, India} \date{}
\maketitle \noindent \doublespacing
\vspace{-0.5cm}
\begin{abstract}
A perfect cuboid, popularly known as a perfect Euler brick/a perfect
box, is a cuboid having integer side lengths, integer face diagonals
and an integer space diagonal. Euler provided an example where only
the body diagonal became deficient for an integer value but it is
known as an Euler brick. Nobody has discovered any perfect cuboid,
however many of us have tried it. The results of this research paper
prove that there exists no perfect cuboid. \\ \it Keywords: {\small
  Perfect Cuboid; Perfect Box; Perfect Euler Brick; Diophantine
  equation.}
\end{abstract}

\section{Introduction}
A cuboid, an Euler brick, is a rectangular parallelepiped with integer
side dimensions together with the face diagonals also as integers. The
earliest time of the problem of finding the rational cuboids can go back
to unknown time, however its existence can be found even before
Euler's work. The definition of an Euler brick in geometric terms can
be formulated mathematically which equivalent to a solution to the
following system of Diophantine equations:
\begin{eqnarray}
a^2+b^2=d^2,~a,b,d\in \mathbb{N};\\
\label{Euler brick 1st}
b^2+c^2=e^2,~b,c,e\in \mathbb{N};\\
\label{Euler brick 2nd}
a^2+c^2=f^2,~a,c,f\in \mathbb{N};
\label{Euler brick 3rd}
\end{eqnarray}
where $a, b, c$ are the edges and $d, e, f$ are the face diagonals. In
1719, Paul Halke \cite{Dickson} found the first known and smallest
Euler Brick with side lengths $\{a,b,c\}=(44,117,240)$ and face
diagonals $\{e,f,g\} = (125,244,267)$. Nicholas Saunderson
\cite{Saunderson}, who was blind from the age of one and the fourth
Lucasian Professor at Cambridge, obtained a parametric solution to the
Euler Brick. Already, families of Euler bricks was announced in 1740
and Euler himself constructed more families of Euler bricks.

Perfect Cuboid is defined as a cuboid where the space diagonal also
has integer length. In other words, the following Diophantine equation
\begin{eqnarray}
a^2+b^2+c^2=g^2;~a,b,c,g\in \mathbb{N}
\label{Perfect Euler brick 1st}
\end{eqnarray}
is included to the system of Diophantine equations (\ref{Euler brick
  1st})-(\ref{Euler brick 3rd}) defining an Euler brick.

There is a question in everyone's mind: ``does a perfect cuboid
exist''? Nobody has discovered any perfect cuboid nor has it been
established that one does not exist, however many of us have tried
it. This problem was in great attention during the 18th century and
Saunderson \cite{Saunderson} reported a parametric solution
\begin{eqnarray}
(a,b,c)=(x(4y^2-z^2),y(x^2-z^2),4xyz)~\text{(if (x,y,z) be a
    Pythagorean tripple)},
\label{Euler brick by Saunderson}
\end{eqnarray}
known as Euler cuboids which always provide Euler bricks although it
does not deliver all possible Euler bricks. We could hope that some of
these Euler cuboids are perfect, but Spohn \cite{Spohn1} demonstrated
that no Euler cuboid can give a perfect cuboid. Althogh, Spohn
\cite{Spohn2} was not completely proved that a derived cuboid of an
Euler cuboid failed to be perfect either, but Chein \cite{Chein} and
Lagrange \cite{Lagrange} both established that this could indeed never
happen. Leech \cite{Leech} provided a one page proof that no Euler
cuboid nor its derived cuboid can be perfect. In 1770 and 1772, Euler
introduced at least two parametric solutions. Euler produced an
example where only the body diagonal falled short of an integer value
(Euler brick). Colman \cite{Colman} shows infinitely many
two-parameter parametrizations of rational cuboids (whose all of the
seven lengths are integers except possibly for one edge (called an
edge cuboid) or face diagonal (called a face cuboid)) with rapidly
increasing degree.

With the help of elementary analysis of the equations for a rational
cuboid modulo for some small primes, Kraitchik \cite{Kraitchik1}
reported that at least one of the sides of a rational cuboid has a
divisor 4 and another one is divided by 16. Moreover, the sides has
divisors as different powers of 3 and at least one of the sides is
divisible by both the primes 5 and 11. The equally elementary
extension of this result was carried out by Horst Bergmann, whereas
Leech showed that the product of all the sides and diagonals (edge and
face) of a perfect cuboid is divisible by $2^8\times 3^4\times
5^3\times 7\times 11\times 13\times 17\times 19\times 29\times 37$
(cf. \cite{Guy}, Problem D18). Kraitchik \cite{Kraitchik1} also
rediscovered the Euler cuboids of (\ref{Euler brick by Saunderson})
and provided a list of 50 rational cuboids which are not Euler cuboids
by using some ad hoc methods. He extended his classical list to 241
cuboids having the odd sides less than $10^6$ in \cite{Kraitchik2} and
found 18 more in \cite{Kraitchik3} out of these 16 were new.

Lal and Blundon \cite{Lal} pointed out that for integers $m,~n,~p$ and
$q$; the cuboid having sides $x=|2mnpq|,~y=|mn(p^2-q^2)|$ and
$z=|pq(m^2-n^2)|$ has at least two face diagonals as integers and is
rational cuboid if and only if $y^2+z^2=\square$. Using the symmetry,
they aimed for a computer search through all the quadruples
$(m,~n,~p,~q)$ satisfying $1\leq m, n, p, q\leq 70$ to check if
$y^2+z^2=\square$ and reported 130 rational cuboids, out of which none
are perfect, however Shanks \cite{Shanks} pointed out some corrigenda
about their paper.

Korec \cite{Korec1} found no perfect cuboids having the least side
smaller than 10000 by the consideration as follows: let $x,~y,~z$ are
the sides of a perfect cuboid, then we can find natural numbers
$a,~b,~c$ all dividing $x$ and $t=\sqrt{y^2+z^2}$ such that
$y=\frac{1}{2}\left(\frac{x^2}{a}-a\right),~z=\frac{1}{2}\left(\frac{x^2}{b}-b\right),~t=\frac{1}{2}\left(\frac{x^2}{c}-c\right)$. Korec
\cite{Korec2} extended his result and found no perfect cuboid if the
least edge smaller than $10^6$. In another research paper, Korec
\cite{Korec3} not found any perfect cuboid if the full diagonal of a
perfect cuboid is $<8\times 10^9$. Moreover, if $x$ and $z$ are the
maximal edge and the full diagonal respectively of a perfect cuboid,
then $z\leq x\sqrt{3}$.

Rathbun \cite{Rathbun1} found 6800 body, 6749 face, 6380 edge and no
perfect cuboids by a computer search if a side $x\leq 333750000$,
while in his another research paper \cite{Rathbun2}, he reported that
4839 of the 6800 body or rational cuboids contain an odd side less
than 333750000 which is a extention and correction of Kraitchik’s
classical table \cite{Kraitchik1,Kraitchik2,Kraitchik3}. If all ``odd
sides''$\leq 10^{10}$, then Butler \cite{Butler} found no perfect cuboids
despite an thorough computer search.

This article is dedicated for the answer of the question ``are there
perfect cuboids?''. It has been discovered that there is no perfect
cuboid.

\section{Results and Discussion}
\subsection{Pythagorean quadruple}
\label{euler_pythagorean_quadruple}
A set of four natural numbers $(a,~b,~c,~d)$ is well known as a
Pythagorean quadruple if the equation $a^2+b^2+c^2=d^2$ satisfies. The
simplest example of a quadruple is $(1,~2,~2,~3)$ as $1^2+2^2+2^2=3^2$
and $(2,~3,~6,~7)$ is the next simplest (primitive) example as
$2^2+3^2+6^2=7^2$. All the primitive quadruples \cite{Spira} can be generated by
the equation
\begin{equation}
(m^2+n^2+p^2+q^2)^2=(m^2+n^2-p^2-q^2)^2+4(mp+nq)^2+4(mq-np)^2;~m,n,p,q\in
  \mathbb{N}.
\label{euler_pythagorean_quadruple_equ}  
\end{equation}
We can find a primitive perfect cuboids iff the following equations are also true.
\begin{eqnarray}
4(mp+nq)^2+4(mq-np)^2=A^2,~A\in \mathbb{N};
\label{euler_brick1}\\
(m^2+n^2-p^2-q^2)^2+4(mp+nq)^2=B^2,~B\in \mathbb{N};
\label{euler_brick2}\\
(m^2+n^2-p^2-q^2)^2+4(mq-np)^2=C^2,~C\in \mathbb{N}.
\label{euler_brick3}
\end{eqnarray}

\subsection{Theorem}
\label{euler_brick_1st_theorem}
For the natural numbers $(a,~b,~c,~d,~e)$; let $(a,~b,~d),~(a,~c,~e)$
be the Pythagorean tripples (three natural numbers)
\begin{eqnarray}
\text{such that}~a^2+b^2=d^2~\text{and}~a^2+c^2=e^2,~\text{then}~~a^2=\frac{(m_2^2-m_1^2)(n_1^2-n_2^2)}{4},~b=\frac{m_1n_1+m_2n_2}{2},~
\label{euler_brick_1st_theorem_equ1}\nonumber\\
c=\frac{m_2n_2-m_1n_1}{2},~d=\frac{m_2n_1+m_1n_2}{2},~e=\frac{m_1n_2-m_2n_1}{2};~\text{where}~m_1,m_2,n_1,n_2\in \mathbb{N}.~~~~~~~~~~~~
\label{euler_brick_1st_theorem_equ2}
\end{eqnarray}

{\bf Proof}: From the equations (\ref{euler_brick_1st_theorem_equ1}), we
can obtain $b^2-c^2=d^2-e^2$. 
\begin{eqnarray}
\text{Or,}~\frac{b-c}{d-e}=\frac{d+e}{b+c}=\frac{m_1}{m_2}
\text{~where $m_1,~m_2~(\in \mathbb{N})$ are relatively prime numbers.}
\label{euler_brick_1st_theorem_equ3}\\
\text{Then}~b-c=\frac{m_1}{m_2}(d-e),~b+c=\frac{m_2}{m_1}(d+e)~~~~~~~~~~~~~~~~~~~~~~~~~~~~~~~~~~~~~~
\label{euler_brick_1st_theorem_equ4}
\end{eqnarray}

Since $b-c$ and $b+c$ are integers and $m_1,~m_2~(\in \mathbb{N})$ are relatively prime numbers, 
\begin{equation}
d-e=m_2n_1~\text{and}~d+e=m_1n_2.~\text{Thus}~d=\frac{m_2n_1+m_1n_2}{2},~e=\frac{m_1n_2-m_2n_1}{2}.
\label{euler_brick_1st_theorem_equ5}
\end{equation}

From equations (\ref{euler_brick_1st_theorem_equ4}) and (\ref{euler_brick_1st_theorem_equ5}), we get
\begin{equation}
b-c=m_1n_1,~b+c=m_2n_2~\text{i.e.}~b=\frac{m_1n_1+m_2n_2}{2},~c=\frac{m_2n_2-m_1n_1}{2}.~~~~~~~~~~~~~~~~~~~~~~~~~~~~~~~~~~~~~~~~~~~~~~~~~~~~~~~~~~~~~~
\label{euler_brick_1st_theorem_equ6}
\end{equation}

Then from equation (\ref{euler_brick_1st_theorem_equ1}),
(\ref{euler_brick_1st_theorem_equ5}) and
(\ref{euler_brick_1st_theorem_equ6}) we get
\begin{equation}
a^2=d^2-b^2=\frac{(m_2^2-m_1^2)(n_1^2-n_2^2)}{4}.~~~~~~~~~~~~~~~~~~~~~~~~~~~~~~~~~~~~~~~~~~~~~~~~~~~~~~~~~~~~~~
\label{euler_brick_1st_theorem_equ7}
\end{equation}

\subsection{Lemma}
\label{euler_brick_1st_lemma}
For the natural numbers $m,n,p,q\in \mathbb{N}$ if
$m^2+n^2=k_1^2r,~p^2+q^2=k_2^2r$ then there exists no perfect cuboid.

{\bf Proof}: From the equation (\ref{euler_brick1}), if
\begin{eqnarray}
4(mp+nq)^2+4(mq-np)^2=4(m^2+n^2)(p^2+q^2)=A^2,~A\in \mathbb{N}
\label{euler_brick_1st_lemma_equ1}\\
\text{with}~m^2+n^2=k_1^2r,~p^2+q^2=k_2^2r\text{ i.e. }m^2p^2+m^2q^2+n^2p^2+n^2q^2=k_1^2k_2^2r^2,
\label{euler_brick_1st_lemma_equ2}\\
\text{then }4(mp+nq)^2+4(mq-np)^2=4k_1^2k_2^2r^2=A^2,~A\in \mathbb{N}.
\label{euler_brick_1st_lemma_equ3}
\end{eqnarray}

{\bf Case (i)}:
\begin{eqnarray}
\text{If}~mp+nq=rk_1k_2 \text{ i.e. }
m^2p^2+n^2q^2+2mnpq=r^2k_1^2k_2^2,~\text{then from equations}
\label{euler_brick_1st_lemma_equ4}~~~~~~~~~~~~~~~~~~~~~~~~~\\
\text{(\ref{euler_brick2}) and
  (\ref{euler_brick_1st_lemma_equ2}), we get }r^2(k_1^2-k_2^2)^2+4r^2k_1^2k_2^2=B^2=r^2(k_1^2+k_2^2)^2,~B\in \mathbb{N}.~~~~~~~~~~~~~~~~~~~~~~~~~~~~~
\label{euler_brick_1st_lemma_equ5}\\
\text{Also from equation (\ref{euler_brick_1st_lemma_equ2}) and
  (\ref{euler_brick_1st_lemma_equ4})},~mq-np=0.~~~~~~~~~~~~~~~~~~~~~~~~~~~~~~~~~~~~~~~
\end{eqnarray}
Thus the equation (\ref{euler_brick3}) has no Pythagorean tripple and
hence, there exists no perfect cuboid for this case.

{\bf Case (ii)}:
\begin{eqnarray}
\text{If}~mq-np=rk_1k_2 \text{ i.e. } m^2q^2+n^2p^2-2mnpq=r^2k_1^2k_2^2\text{ then from equations}
\label{euler_brick_1st_lemma_equ6}~~~~~~~~~~~~~~~~~~~~~~~~~\\
\text{(\ref{euler_brick_1st_lemma_equ2}) and (\ref{euler_brick3}), we
  get } r^2(k_1^2-k_2^2)^2+4r^2k_1^2k_2^2=C^2=r^2(k_1^2+k_2^2)^2,~C\in
\mathbb{N}.~~~~~~~~~~~~~~~~~~~~~~~~~~~~~
\label{euler_brick_1st_lemma_equ7}\\
\text{Also from equation (\ref{euler_brick_1st_lemma_equ2}) and
  (\ref{euler_brick_1st_lemma_equ6})},~mp+nq=0.~~~~~~~~~~~~~~~~~~~~~~~~~~~~~~~~~~~~~~~
\end{eqnarray}
Hence the equation (\ref{euler_brick2}) has no Pythagorean tripple and
hence, there exists no perfect cuboid for this case.

{\bf Case (iii)}:

Let us consider $m_1=k_1$, $m_2=k_2$, $n_1=2rk_2$ and $n_2=2rk_1$ in
the equations (\ref{euler_brick_1st_theorem_equ5}),
(\ref{euler_brick_1st_theorem_equ6}) and
(\ref{euler_brick_1st_theorem_equ7}); then $a^2=r^2(k_2^2-k_1^2)^2$,
$b=2rk_1k_2$, $c=0$. Then the equation (\ref{euler_brick3}) has no
Pythagorean tripple and hence, there exists no perfect cuboid for this
case.

\subsection{Lemma}
\label{euler_brick_2nd_lemma}
For the natural numbers $m,n,p,q\in \mathbb{N}$ if
$m^2+n^2=k_1^2,~p^2+q^2=k_2^2$; then there exists no perfect cuboid.

{\bf Proof}: From the equation (\ref{euler_brick1}), if
\begin{eqnarray}
4(mp+nq)^2+4(mq-np)^2=4(m^2+n^2)(p^2+q^2)=A^2,~A\in \mathbb{N}
\label{euler_brick_2nd_lemma_equ1}\\
\text{with}~m^2+n^2=k_1^2,~p^2+q^2=k_2^2\text{ i.e. }m^2p^2+m^2q^2+n^2p^2+n^2q^2=k_1^2k_2^2,
\label{euler_brick_2nd_lemma_equ2}\\
\text{then }4(mp+nq)^2+4(mq-np)^2=4k_1^2k_2^2=A^2,~A\in \mathbb{N}
\label{euler_brick_2nd_lemma_equ3}
\end{eqnarray}

{\bf Case (i)}:
\begin{eqnarray}
\text{If}~mp+nq=k_1k_2 \text{ i.e. }
m^2p^2+n^2q^2+2mnpq=k_1^2k_2^2,~\text{then from equations}
\label{euler_brick_2nd_lemma_equ4}~~~~~~~~~~~~~~~~~~~~~~~~~\\
\text{(\ref{euler_brick2}) and
  (\ref{euler_brick_2nd_lemma_equ2}), we get }(k_1^2-k_2^2)^2+4k_1^2k_2^2=B^2=(k_1^2+k_2^2)^2,~B\in \mathbb{N}.~~~~~~~~~~~~~~~~~~~~~~~~~~~~~
\label{euler_brick_2nd_lemma_equ5}\\
\text{Also from equation (\ref{euler_brick_2nd_lemma_equ2}) and
  (\ref{euler_brick_2nd_lemma_equ4})},~mq-np=0.~~~~~~~~~~~~~~~~~~~~~~~~~~~~~~~~~~~~~~~
\end{eqnarray}
Thus the equation (\ref{euler_brick3}) has no Pythagorean tripple and
hence, there exists no perfect cuboid for this case.

{\bf Case (ii)}:
\begin{eqnarray}
\text{If}~mq-np=k_1k_2 \text{ i.e. } m^2q^2+n^2p^2-2mnpq=k_1^2k_2^2\text{ then from equations}
\label{euler_brick_2nd_lemma_equ6}~~~~~~~~~~~~~~~~~~~~~~~~~\\
\text{(\ref{euler_brick_2nd_lemma_equ2}) and (\ref{euler_brick3}), we
  get } (k_1^2-k_2^2)^2+4k_1^2k_2^2=C^2=(k_1^2+k_2^2)^2,~C\in
\mathbb{N}.~~~~~~~~~~~~~~~~~~~~~~~~~~~~~
\label{euler_brick_2nd_lemma_equ7}\\
\text{Also from equation (\ref{euler_brick_2nd_lemma_equ2}) and
  (\ref{euler_brick_2nd_lemma_equ6})},~mp+nq=0.~~~~~~~~~~~~~~~~~~~~~~~~~~~~~~~~~~~~~~~
\end{eqnarray}
Hence the equation (\ref{euler_brick2}) has no Pythagorean tripple and
hence, there exists no perfect cuboid for this case.

{\bf Case (iii)}:

Let us consider $m_1=k_1$, $m_2=k_2$, $n_1=2k_2$
and $n_2=2k_1$ in the equations (\ref{euler_brick_1st_theorem_equ5}),
(\ref{euler_brick_1st_theorem_equ6}) and
(\ref{euler_brick_1st_theorem_equ7}); then $a^2=(k_2^2-k_1^2)^2$,
$b=2k_1k_2$, $c=0$. Then the equation (\ref{euler_brick3}) has no
Pythagorean tripple and hence, there exists no perfect cuboid for this
case.

\subsection{Lemma}
\label{euler_brick_3rd_lemma}
For the natural numbers $m,n,p,q\in \mathbb{N}$ if
$m^2+n^2=r,~p^2+q^2=k_2^2r$ then there exists no perfect cuboid.

{\bf Proof}: From the equation (\ref{euler_brick1}), if
\begin{eqnarray}
4(mp+nq)^2+4(mq-np)^2=4(m^2+n^2)(p^2+q^2)=A^2~(A\in \mathbb{N}),
\label{euler_brick_3rd_lemma_equ1}\\
\text{with}~m^2+n^2=r,~p^2+q^2=k_2^2r\text{ i.e. }m^2p^2+m^2q^2+n^2p^2+n^2q^2=r^2k_2^2,
\label{euler_brick_3rd_lemma_equ2}\\
\text{then }4(mp+nq)^2+4(mq-np)^2=4k_2^2r^2=A^2,~A\in \mathbb{N}
\label{euler_brick_3rd_lemma_equ3}
\end{eqnarray}

{\bf Case (i)}:
\begin{eqnarray}
\text{If}~mp+nq=rk_2 \text{ i.e. }
m^2p^2+n^2q^2+2mnpq=r^2k_2^2,~\text{then from equations}
\label{euler_brick_3rd_lemma_equ4}~~~~~~~~~~~~~~~~~~~~~~~~~\\
\text{(\ref{euler_brick2}) and (\ref{euler_brick_3rd_lemma_equ2}), we
  get }r^2(1-k_2^2)^2+4k_1^2k_2^2=B^2=r^2(1+k_2^2)^2,~B\in
\mathbb{N}.~~~~~~~~~~~~~~~~~~~~~~~~~~~~~
\label{euler_brick_3rd_lemma_equ5}\\
\text{Also from equation (\ref{euler_brick_3rd_lemma_equ2}) and
  (\ref{euler_brick_3rd_lemma_equ4})},~mq-np=0.~~~~~~~~~~~~~~~~~~~~~~~~~~~~~~~~~~~~~~~
\end{eqnarray}
Thus the equation (\ref{euler_brick3}) has no Pythagorean tripple and
hence, there exists no perfect cuboid for this case.

{\bf Case (ii)}:
\begin{eqnarray}
\text{If}~mq-np=rk_2 \text{ i.e. } m^2q^2+n^2p^2-2mnpq=r^2k_2^2\text{ then from equations}
\label{euler_brick_3rd_lemma_equ6}~~~~~~~~~~~~~~~~~~~~~~~~~\\
\text{(\ref{euler_brick_3rd_lemma_equ2}) and (\ref{euler_brick3}), we
  get } r^2(1-k_2^2)^2+4r^2k_2^2=C^2=r^2(1+k_2^2)^2,~C\in
\mathbb{N}.~~~~~~~~~~~~~~~~~~~~~~~~~~~~~
\label{euler_brick_3rd_lemma_equ7}\\
\text{Also from equation (\ref{euler_brick_3rd_lemma_equ2}) and
  (\ref{euler_brick_3rd_lemma_equ6})},~mp+nq=0.~~~~~~~~~~~~~~~~~~~~~~~~~~~~~~~~~~~~~~~
\end{eqnarray}
Hence the equation (\ref{euler_brick2}) has no Pythagorean tripple and
hence, there exists no perfect cuboid for this case.

{\bf Case (iii)}:

Let us consider $m_1=k_2$, $m_2=1$, $n_1=2r$ and
$n_2=2rk_2$ in the equations (\ref{euler_brick_1st_theorem_equ5}),
(\ref{euler_brick_1st_theorem_equ6}) and
(\ref{euler_brick_1st_theorem_equ7}); then $a^2=r^2(1-k_2^2)^2$,
$b=2rk_2$, $c=0$. Then the equation (\ref{euler_brick3}) has no
Pythagorean tripple and hence, there exists no perfect cuboid for this
case.

\subsection{Lemma}
\label{euler_brick_4th_lemma}
For the natural numbers $m,n,p,q\in \mathbb{N}$ if
$m^2+n^2=k_1^2r,~p^2+q^2=r$; then there exists no perfect cuboid.

{\bf Proof}: From the equation (\ref{euler_brick1}), if
\begin{eqnarray}
4(mp+nq)^2+4(mq-np)^2=4(m^2+n^2)(p^2+q^2)=A^2,~A\in \mathbb{N}
\label{euler_brick_4th_lemma_equ1}\\
\text{with}~m^2+n^2=k_1^2r,~p^2+q^2=r\text{ i.e. }m^2p^2+m^2q^2+n^2p^2+n^2q^2=k_1^2r^2,
\label{euler_brick_4th_lemma_equ2}\\
\text{then }4(mp+nq)^2+4(mq-np)^2=4k_1^2r^2=A^2,~A\in \mathbb{N}
\label{euler_brick_4th_lemma_equ3}
\end{eqnarray}

{\bf Case (i)}:
\begin{eqnarray}
\text{If}~mp+nq=rk_1 \text{ i.e. }
m^2p^2+n^2q^2+2mnpq=r^2k_1^2,~\text{then from equations}
\label{euler_brick_4th_lemma_equ4}~~~~~~~~~~~~~~~~~~~~~~~~~\\
\text{(\ref{euler_brick2}) and
  (\ref{euler_brick_4th_lemma_equ2}), we get }r^2(k_1^2-1)^2+4r^2k_1^2=B^2=r^2(k_1^2+1)^2,~B\in \mathbb{N}.~~~~~~~~~~~~~~~~~~~~~~~~~~~~~
\label{euler_brick_4th_lemma_equ5}\\
\text{Also from equation (\ref{euler_brick_4th_lemma_equ2}) and
  (\ref{euler_brick_4th_lemma_equ4})},~mq-np=0.~~~~~~~~~~~~~~~~~~~~~~~~~~~~~~~~~~~~~~~
\end{eqnarray}
Thus the equation (\ref{euler_brick3}) has no Pythagorean tripple and
hence, there exists no perfect cuboid for this case.

{\bf Case (ii)}:
\begin{eqnarray}
\text{If}~mq-np=rk_1 \text{ i.e. } m^2q^2+n^2p^2-2mnpq=r^2k_1^2\text{ then from equations}
\label{euler_brick_4th_lemma_equ6}~~~~~~~~~~~~~~~~~~~~~~~~~\\
\text{(\ref{euler_brick_4th_lemma_equ2}) and (\ref{euler_brick3}), we
  get } r^2(k_1^2-1)^2+4r^2k_1^2=C^2=r^2(k_1^2+1)^2,~C\in
\mathbb{N}.~~~~~~~~~~~~~~~~~~~~~~~~~~~~~
\label{euler_brick_4th_lemma_equ7}\\
\text{Also from equation (\ref{euler_brick_1st_lemma_equ2}) and
  (\ref{euler_brick_4th_lemma_equ6})},~mp+nq=0.~~~~~~~~~~~~~~~~~~~~~~~~~~~~~~~~~~~~~~~
\end{eqnarray}
Hence the equation (\ref{euler_brick2}) has no Pythagorean tripple and
hence, there exists no perfect cuboid for this case.

{\bf Case (iii)}:

Let us consider $m_1=1$, $m_2=k_1$, $n_1=2rk_1$ and
$n_2=2r$ in the equations (\ref{euler_brick_1st_theorem_equ5}),
(\ref{euler_brick_1st_theorem_equ6}) and
(\ref{euler_brick_1st_theorem_equ7}); then $a^2=r^2(k_1^2-1)^2$, $b=2rk_1$,
$c=0$. Then the equation (\ref{euler_brick3}) has no Pythagorean
tripple and hence, there exists no perfect cuboid for this case.

\subsection{Remark}
\label{euler_brick_1st_remark}
 The results of other cases for the values of the $m_1$, $m_2$,
 $n_1$, $n_2$ regarding the Lemmas
 \ref{euler_brick_1st_lemma}-\ref{euler_brick_4th_lemma} are straight
 forward. For example, if you consider the values of the variables as
 $m_1=rk_1$, $m_2=rk_2$, $n_1=2k_2$, $n_1=2k_1$ in the equations
 (\ref{euler_brick_1st_theorem_equ5}),
 (\ref{euler_brick_1st_theorem_equ6}) and
 (\ref{euler_brick_1st_theorem_equ7}); then $a^2=r^2(k_2^2-k_1^2)^2$,
 $b=2rk_1k_2$, $c=0$. Then the equation (\ref{euler_brick3}) has no
 Pythagorean tripple and hence, there exists no perfect cuboid for this
 case.

\subsection{Lemma}
\label{euler_brick_5th_lemma}
For the natural numbers $m,n,p,q\in \mathbb{N}$,
$m^2+n^2=p^2+q^2$, then there exists no perfect cuboid.

{\bf Proof}: If the natural numbers $m,n,p,q\in \mathbb{N}$ such that
$m^2+n^2=p^2+q^2$, then from equations
(\ref{euler_pythagorean_quadruple_equ}), (\ref{euler_brick2}) and
(\ref{euler_brick3}) we can't get any Pythagorean quadruple and
Pythagorean tripple respectively and hence, there exists no perfect cuboid for this case.

\subsection{Theorem}
\label{euler_brick_2nd_theorem}
There exists no perfect cuboid.

{\bf Proof}: We can discover any perfect cuboid only when we can get
any natural number solutions of equations
(\ref{euler_pythagorean_quadruple_equ})-(\ref{euler_brick3}). Hence by
the the Lemmas
\ref{euler_brick_1st_lemma}-\ref{euler_brick_5th_lemma}, we can
conclude that there exists no perfect cuboid.

\vspace{1cm}

{\bf Acknowledgment:} {\it I am grateful to the University Grants
  Commission (UGC), New Delhi for awarding the Dr. D. S. Kothari Post
  Doctoral Fellowship from 9th July, 2012 to 8th July, 2015 at Indian
  Institute of Technology (BHU), Varanasi. The self training for this
  investigation was continued during the period.}

\end{document}